\documentclass[smallextended]{svjour3}       % onecolumn (second format)

\smartqed  % flush right qed marks, e.g. at end of proof

\usepackage{graphicx}

\usepackage{mathptmx}      % use Times fonts if available on your TeX system

\usepackage{amsmath}
\usepackage{amssymb}
\usepackage{amsthm}
\usepackage[nodate]{datetime}
\usepackage{dsfont}
\usepackage{fancyhdr}
\usepackage[latin1]{inputenc}
\usepackage{psfrag}
\usepackage{srcltx}
\usepackage{stmaryrd}
\usepackage[normalem]{ulem}
\usepackage{url}

\newcommand{\PPP}{\mathcal P}
\newcommand{\N}{\mathbb{N}}
\newcommand{\Z}{\mathbb{Z}}

\newcommand{\Pb}{\mathbb{P}}
\newcommand{\I}{\mathds{1}}
\newcommand{\II}{\mathcal I}

\newcommand{\h}{h}

\psfrag{x1}{\small $x_1$}
\psfrag{a1}{\small $a_1$}
\psfrag{x2}{\small $x_2$}
\psfrag{a2}{\small $a_2$}
\psfrag{?} {\small ?}
\psfrag{0} {\small $0$}

\theoremstyle{plain}
\newtheorem{theo}{Theorem}
\newtheorem*{theo*}{Theorem}

\newtheorem{lemma}{Lemma}
\newtheorem{cor}[lemma]{Corollary}

\theoremstyle{definition}
\theoremstyle{remark}

\newtheorem*{claim*}{Claim}

\newtheorem*{remark*}{Remark}
\newtheorem{property}{Property}

\journalname{}

\begin{document}

\title{Absorbing-State Phase Transition for Driven-Dissipative Stochastic
Dynamics on $\mathbb Z$}

\author{Leonardo T.\ Rolla and Vladas Sidoravicius \\
{\small ENS-Paris and IMPA}}

% \subtitle{Do you have a subtitle?\\ If so, write it here}

\titlerunning{Stochastic Sandpiles and Activated Random Walks}        % if too long for running head

\authorrunning{Leonardo T. Rolla and Vladas Sidoravicius} % if too long for running head

\institute{L. T. Rolla \at
              École Normale Supérieure.
              45, rue d'Ulm - Paris 75005 - France\\
              Instituto de Matemática Pura e Aplicada.
              Estrada Dona Castorina, 110 - Rio de Janeiro 22.460-320 - Brasil
           \and
              V. Sidoravicius \at
              Instituto de Matemática Pura e Aplicada.
              Estrada Dona Castorina, 110 - Rio de Janeiro 22.460-320 - Brasil
}

\date{August 30, 2011}

\maketitle

\begin{abstract}
We study the long-time behavior of conservative interacting particle systems in $\mathbb Z$:
the activated random walk model for reaction-diffusion systems and the stochastic sandpile.
We prove that both systems undergo an absorbing-state phase transition.
\\

\noindent
This preprint has the same numbering for sections, theorems, equations and figures as the published article \emph{``Invent. Math. 188 (2012): 127--150.''}
\\
\end{abstract}

\section{Introduction}
\label{sec:sec1introduction}

Modern Statistical Mechanics offers a large and important class of driven-dissipative lattice systems that naturally evolve to a critical state, which is characterized by power-law distributions of the sizes of relaxation events (a paradigm example is the emergence of avalanches caused by small perturbations).
In many mathematically interesting and physically relevant cases such systems are attracted to a stationary critical state without being specifically tuned to a critical point.
In particular, it is believed that this phenomenon lies behind random fluctuations at the macroscopic scale, and creation of self-similar shapes in a variety of growth systems.

Due to strong non-locality of correlations and dynamic long-range effects, classical analytic and probabilistic techniques fail in most cases of interest, making the rigorous analysis of such systems a major mathematical challenge.

Studies of the above phenomenon are confined to very few models, and its conceptual understanding is extremely fragmented.
Among theories which attempt to explain long-ranged spatio-temporal correlations, the physical paradigm called `self-organized criticality' takes its particular place~\cite{dhar-06,hinrichsen-00}.
These are systems whose natural dynamics drives them towards, and then maintains them at the \emph{edge} of stability~\cite{dhar-06}.
However, for non-equilibrium steady states it is becoming increasingly evident that `self-organized criticality' is related to conventional critical behavior, namely that of an \emph{absorbing-state phase transition}.
The known examples are variations of underlying non-equilibrium systems which actually do have a parameter and exhibit critical phenomena.
The phase transition in these systems arises from a conflict between a spread of activity and a tendency for this activity to die out~\cite{dickman-02,lubeck-04}, and the transition point separates an active and an absorbing phase in which the dynamics gets eventually extinct in any finite region.
\medskip

In this paper we focus on two chief examples of conservative, infinite-volume systems which belong to the above mentioned family: the activated random walk model for reaction-diffusion and the stochastic sandpile model.
Let us briefly describe the models (see Section~\ref{sec2results} for the precise definitions).

The reaction-diffusion model is given by the following conservative particle dynamics in $\Z^d$.
Each particle can be in one of two states: an active $A$-state, and a passive $S$-state.
Particles in the $A$-state perform a continuous-time random walk with jump rate $D_A=1$ without interacting.
Particles in the $S$-state do not move, that is, $D_S=0$.
Each particle changes its state $A \to S$ at some halting rate $\lambda>0$ and the reaction $A+S \to 2A$ happens immediately.
The catalyzed transition $A+S \to 2A$ and the spontaneous transition $A\to S$ represent the spread of activity versus a tendency for this activity to die out.
This system will be referred to as the model of Activated Random Walks (ARW).\footnote{The Activated Random Walks may be viewed as a special case of a {diffusive epidemic process}.
The model was introduced in the late 1970's by F.~Spitzer, but due to its tremendous technical difficulties and complexity, remained unsolved until recently, when it was studied in detail in~\cite{kesten-sidoravicius-03,kesten-sidoravicius-05,kesten-sidoravicius-06,kesten-sidoravicius-08}.
This process has also been studied via renormalization group techniques and numerical simulation~\cite{defreitas-lucena-dasilva-hilhorst-00,fulco-messias-lyra-01,fulco-messias-lyra-01-1,janssen-01,kree-schaub-schmittmann-89,oerding-vanwijland-leroy-hilhorst-00,vanwijland-oerding-hilhorst-98}.}

In sandpile models the state of the system is represented by the number of particles $\eta(x)=0,1,2,\dots$ at each vertex $x\in\Z^d$.
The vertex $x$ is \emph{stable} when $\eta(x)<N_c$, for some threshold value $N_c$, and \emph{unstable} when $\eta(x) \geqslant N_c$.
Relaxation (update of the state) happens by toppling each unstable vertex, i.e., sending particles to its neighbors following a certain (deterministic or stochastic) rule.
In this paper we study the following sandpile dynamics, which is a variation of Manna's model frequently considered in the physics literature~\cite{dhar-99,manna-91}.
The threshold for stability of vertices is $N_c=2$, and each unstable vertex topples after an exponentially-distributed time, sending $2$ particles to neighbors chosen independently at random.
We will refer to this model as the Stochastic Sandpile Model (SSM).\footnote{In the so-called \emph{Abelian Sandpile}, $N_c = 2d$, and an unstable vertex \emph{deterministically} sends one particle to each of its $2d$ neighbors when toppling~\cite{bak-tang-wiesenfeld-87,bak-tang-wiesenfeld-88}.}
\medskip

For these systems, the relation between self-organized and ordinary criticality is understood as follows.
On the one hand, `self-organized criticality' appears in their parameter-free,{\footnotemark}\ finite-volume variation: particles are added to the bulk of a finite box, and absorbed at its boundary during relaxation.
The particle addition happens at a \emph{slow} rate, or \emph{only after} the system globally stabilizes.
In this dynamics, when the average density $\mu$ inside the box is too small, mass tends to accumulate.
When it is too large, there is intense activity and a substantial number of particles is absorbed at the boundary.
With this \emph{carefully designed} mechanism, the model is attracted to a \emph{critical state} with an average density given by $0<\mu_c<\infty$, though it was not explicitly tuned to this critical value.
On the other hand, the corresponding conservative systems in infinite volume exhibit ordinary criticality in the sense that their dynamics fixate for $\mu<\mu_c$ and do not fixate for $\mu>\mu_c$, and moreover the \emph{critical exponents} of the finite-volume addition-relaxation dynamics are related to those of the conservative dynamics in infinite volume.

\footnotetext{Consider the \emph{density of particles} as a parameter, and let $\lambda$ be just a fixed number.}

The \emph{critical behavior} of stochastic sandpiles seems to belong to the same \emph{universality class} as the depinning of a linear elastic interface subject to random pinning potentials,\footnote
{This prediction is supported by simulations, at least for dimensions $d\geqslant2$, see~\cite{dickman-02}.}
roughly depicted by a nailed carpet being detached from the floor by an external force of critical intensity, where the rupture of each nail induces other ruptures nearby, giving rise to ``avalanches''.

The deterministic sandpile defines a universality class \textsl{sui generis}, and is marked by strong non-ergodic effects~\cite{dickman-02}.\footnote
{\label{asm-note}%
For instance, it is known~\cite{feydenboer-redig-05,meester-quant-05,redig-06} that there is local fixation for $\mu<d$ and there is no fixation for $\mu>2d-1$, but one can construct somewhat artificial ergodic states with $\mu\in(\mu_c,2d-1)$ that fixate and there are as well states with $\mu\in(d,\mu_c)$ that do not fixate~\cite{feydenboer-redig-05}.
Furthermore, as observed in~\cite{fey-levine-wilson-10-1,fey-levine-wilson-10}, the density $\mu_s$ attained in the large-time limit for the driven-dissipative system does not even satisfy $\mu_s=\mu_c$.
}
This seems to be due to the failure of the toppling procedure in eliminating certain microscopic symmetries in the configuration by the time the system becomes unstable, as a consequence of the existence of many toppling invariants.
At the same time, the stationary state of the deterministic sandpile is the uniform distribution over a special subset of configurations exhibiting several combinatorial properties which, together with its strong algebraic structure, is one of the reasons why most of the mathematical literature has focused on this model (see~\cite{redig-06} and references therein).

In the stochastic sandpile models, by contrast, the addition-relaxation operators are themselves random, leading to a set of coupled polynomial equations~\cite{dhar-06,sadhu-dhar-09}.
Their explicit solutions are not known in a general form, and very little can be said rigorously about the phase transition of such systems.
The reaction-diffusion dynamics, in turn, might be even more apt to dispel microscopic symmetries, for not only the moves are random, but also the particles jump individually rather than in pairs.
\medskip

The main pursuit in this framework is to describe the critical behavior, the scaling relations and critical exponents, and whether the critical density is the same as the long-time limit attained in the driven-dissipative version.
These questions are however far beyond the reach of current techniques.

The contribution of this paper is a step forward in the mathematical understanding of these systems: we develop a general method for the analysis of diffusive-dissipative dynamics, and establish existence of the absorbing-state phase transition for both models in the one-dimensional case.
Namely, we show that the system fixates if it starts with small enough density.
Remark that this question remains open for $d \geqslant 2$ in spite of considerable efforts in the probability and mathematical physics communities.

\paragraph{A word about the proof.}
We build upon the Diaconis-Fulton representation to exploit the combinatorial nature of the problem.
Due to particle exchangeability, this representation extracts precisely the part of the randomness that is relevant for the phase transition, focusing on the total number of jumps and leaving aside the order in which they take place.
It is suitable for studying path traces, total occupation times, and final particle positions, but precludes the analysis of quantities for which the order of the jumps does matter, such as correlation functions or local shape properties.

We introduce a tool to modify the
Diaconis-Fulton instructions, which consists in forcing
stable particles to move whereas in the original dynamics they would stand
still, and we show that the total amount of activity is non-decreasing with
respect to this change.
For the ARW model this means suppressing some $A\to S$ transitions and for the SSM it is done via extra ``semi-legal'' single-particle topplings.

The instruction manipulations give broad possibilities of handling the
configuration as it evolves, for one can keep the particles moving until they
arrive at more convenient locations.
We achieve fixation by successively settling the particles at appropriate
places.
In this procedure we set up a trap for each particle,
yielding a globally stable configuration.
We finally show that this scheme has positive probability of working indefinitely.

An important step is to look ahead in the future before deciding which semi-legal operations should be performed.
It is crucial though that the Diaconis-Fulton instructions are explored in a careful way, in order to keep independence.
This approach allows a fine control over the motion of each individual particle.
\medskip

This paper is organized as follows.
In Section~\ref{sec2results} we give precise definitions of the models and statements of the results.
In Section~\ref{sec1flag} we study the Diaconis-Fulton representation and relate stability of the initial configuration to local fixation of the system.
In Section~\ref{sec0heuristics} we explain the general strategy used in the proofs and give an overview of the main ideas.
In Sections~\ref{sec1phasetransition} and~\ref{sec5phasetransitionssm} we prove the phase transition respectively for the ARW and SSM.
In Section~\ref{sec1concluding} we discuss extensions of our results and conclude with open problems.

\section{The models and results}
\label{sec2results}

\paragraph{Stochastic Sandpile Model.}
When a site $x$ has at least $2$ particles, it is called \emph{unstable} and topples at rate $1$.
When it topples it sends $2$ particles to neighboring sites chosen independently at random, that is, according to the distribution $p(y-x)$, where $p(z)=\frac1{2d}$ if $\|z\|=1$ and $0$ otherwise.

The state of the SSM at each time $t\geqslant0$ is given by $\eta_t\in(\N_{0})^{\Z^d}$, where $\N_0 = \N \cup \{0\}$ and $\eta_t(x)$ denotes the number of particles found at site $x$ at time $t$.
For each site $x\in\Z^d$, the transitions $\eta\to\tau_{xy}\tau_{xw}\eta$ happen at rate $A\big(\eta_t(x)\big)p(y-x)p(w-x)$, where
\[
  \tau_{xy}\eta(z) =
  \begin{cases}
    \eta(x)-1, & z=x \\
    \eta(y)+1, & z=y \\
    \eta(z),   & \mbox{otherwise,}
  \end{cases}
\]
and $A(k)=\I_{k\geqslant 2}$ indicates whether or not site $x$ is unstable.
Let $\Pb^\nu$ denote the law of $(\eta_t)_{t\geqslant0}$ starting from $\eta_0$ distributed as $\nu$.
This evolution is well defined because the jump rates are bounded.

We say that the system \emph{locally fixates} if $\eta_t(x)$ is eventually constant for each $x$, otherwise we say that the system \emph{stays active}.

\begin{theo}
\label{theo1fixationssm}
Consider the Stochastic Sandpile Model in the one-dimensional lattice $\Z$, with initial distribution $\nu$ given by i.i.d.\ Poisson random variables with parameter $\mu$.
There exists $\mu_c \in \big[\frac14 ,1\big]$ such that the system locally fixates a.s.\ if $\mu < \mu_c$, and stays active a.s.\ if $\mu > \mu_c$.
\end{theo}

\paragraph{Activated Random Walk model.}
Each particle in the $A$-state performs a continuous-time random walk with jump rate $D_A=1$.
The jumps have a probability density $p(\cdot)$ on $\Z^d$ such that the set $\{z\in\Z^d:p(z)>0\}$ generates the whole group $(\Z^d,+)$.
Independently of anything else, each particle in the $A$-state turns to the $S$-state at a halting rate $\lambda>0$.
Once a particle is in the $S$-state, it stops moving, i.e., its jump rate is $D_S=0$, and it remains in the $S$-state until the instant when another particle is present at the same vertex.
At such an instant the particle which is in $S$-state flips to the $A$-state, giving the transition $A+S \to 2A$.
A particle in the $S$-state stands still forever if no other particle ever visits the vertex where it is located.
According to these rules, the transition $A \to S$ {\it effectively} occurs if and only if, at the instant of such a transition, the particle \emph{does not share} its vertex with another particle (the innocuous instantaneous transition $2A \to A+S \to 2A$ is not observed).
Particles in the $A$-state do not interact among themselves.

The state of the ARW at time $t\geqslant 0$ is given by $\eta_t \in \Sigma = (\N_{0\varrho})^{\Z^d}$, where $\N_{0\varrho} = \N_0 \cup \{\varrho\}$.
In this setting $\eta_t(x)$ denotes the \emph{number of particles} at site $x$ at time $t$, i.e., the \emph{number of active particles} if $\eta_t(x)\in\N$, and \emph{one passive particle} if $\eta_t(x)=\varrho$.
We turn $\N_{0\varrho}$ into an ordered set by setting $0<\varrho<1<2<\cdots$.
We also let $|\varrho|=1$, so $|\eta_t(x)|$ counts the number of particles regardless of their state.
The addition is defined by $\varrho+0=0+\varrho=\varrho$ and $\varrho+n=n+\varrho=1+n$ for $n>0$, that is, the addition operation already provides the $A+S\to2A$ transition.
The $A\to S$ transition is represented by $\varrho\cdot n$ given by $\varrho\cdot 1=\varrho$ and $\varrho\cdot n=n$ for $n>1$.
Subtractions involving $\varrho$ are not defined.

The process evolves as follows.
For each site $x$, we have the transitions $\eta\to \tau_{xy}\eta$ at rate $A\big(\eta_t(x)\big)p(y-x)$ and $\eta\to \tau_{x\varrho}\eta$ at rate $\lambda A\big(\eta_t(x)\big)$, where
\[
  \tau_{x\varrho}\eta(z) =
  \begin{cases}
    \varrho\cdot\eta(x), & z=x \\
    \eta(z),   & \mbox{otherwise}
  \end{cases}
\]
and $A(k)=k\I_{k\geqslant1}$, so $A\big(\eta_t(x)\big)$ is the number of active particles at site $x$ at time $t$.
We denote by $\nu$ the distribution of $\eta_0$ and assume that $\eta_0(x)\in\N_0$ for all $x$ a.s.
We further write $\nu_M$ for the distribution of the truncated configuration $\eta^M$ given by $\eta^M(x)=\eta_0(x)$ for ${|x|\leqslant M}$ and $\eta^M(x)=0$ otherwise, and $\Pb^\nu_M=\Pb^{\nu_M}$.
$\Pb^\nu_M$ is well defined and corresponds to the evolution of a
countable-state Markov chain whose configurations contain only finitely many
particles.

It follows from a construction due to Andjel\footnote
  {Let $\Sigma'' = \{\eta\in\Sigma: \sum_x |\eta(x)| <\infty \}$,
  define $\alpha(x)=\sum_{n=0}^\infty 2^{-n}p^{(n)}(x,0)$
  and
  $\|\eta\| = \sum_{x\in \Z^d} |\eta(x)|\alpha(x)$.
  Let $\Sigma'=\{\eta\in\Sigma:\|\eta\|<\infty\}$.
  Then $\Sigma''$ is dense in $\Sigma'$, $\nu(\Sigma')=1$,
  and a standard adaptation of E.~Andjel's construction (see~\cite{andjel-82})
  implies the existence of a measure $\Pb^\nu$ on $\Omega=D\bigl([0,\infty),\Sigma')$
  corresponding to a Markov process with the prescribed transition rates,
  having the Feller property in this topology, and thus satisfying~\eqref{eq1finitebox}.}
that, if $\nu$ is a product measure with density $\nu(|\eta(0)|)<\infty$ then $\Pb^\nu$ is well defined and, moreover,
\begin{equation}
\label{eq1finitebox}
  \Pb^\nu(E) = \lim_{M\to\infty} \Pb^\nu_{M}(E)
\end{equation}
for every event $E$ that depends on a finite space-time window, i.e., that is measurable with respect to $\big(\eta_s(x):{|x|<t,s\in[0,t]}\big)$ for some $t<\infty$.

\begin{theo}
\label{theo1fixation}
Consider the Activated Random Walk Model with nearest-neighbor jumps in the one-dimensional lattice $\Z$ with halting rate $\lambda$ and with initial distribution $\nu$ given by i.i.d.\ Poisson random variables in $\N_0$ with parameter $\mu$.
There exists $\mu_c \in \big[\frac\lambda{1+\lambda} ,1\big]$ such that the system locally fixates a.s.\ if $\mu < \mu_c$ and stays active a.s.\ if $\mu > \mu_c$.
\end{theo}

Theorem~\ref{theo1fixation} should be contrasted with the particular case of totally asymmetric jumps, for which it is known~\cite{hoffman-sidoravicius-04} that $\mu_c=\frac{\lambda}{1+\lambda}$, so the lower bound is sharp.
Note also that $\mu_c = 1$ when $\lambda = + \infty$.
Theorems~\ref{theo1fixationssm}~and~\ref{theo1fixation} remain true under more general hypotheses.
See Section~\ref{sec1concluding} for details.

\section{The Diaconis-Fulton representation}
\label{sec1flag}

In this section we describe the Diaconis-Fulton graphical representation for the dynamics of both SSM and ARW.
The main advantage of this representation is that it leaves aside the chronological order of topplings.
This will be justified \textsl{a posteriori} by showing that local fixation for the stochastic evolution is equivalent to stability properties of the discrete operations described hereafter.

Graphical constructions have become standard tools in Probability.
They are usually attributed to Harris, see~\cite{harris-78} and references therein to backtrack its history.
In the construction described below, a sequence of random instructions is
assigned to each site, but unlike Harris' construction, the instant when each
step should be performed is not part of the representation and is determined by
some external factor.
As far as we know, the earliest references describing this kind of construction are the card game of Diaconis and Fulton~\cite{diaconis-fulton-91} and the general framework considered by Eriksson~\cite{eriksson-96}, which includes the card game, the Abelian Sandpile Model, the Stochastic Sandpile Model, and the Activated Random Walks.

We first describe the \emph{representation}, defining $\II$, $\PPP$, $\h$, $\Phi$, instructions, stable sites and legal operations.
We then stress some \emph{local} properties in an \emph{abstract}, model-independent framework, defining $m_{\alpha}$, $\Phi_\alpha$, and legal sequences.
We finally consider the \emph{global} consequences of these properties that will
be central in the upcoming sections, and define $m_{V,\eta}$, stable
configurations and stabilizing sequences.

\paragraph{Representation.}
We start with the representation for the SSM.
A site $x\in\Z^d$ is \emph{stable} in the configuration $\eta$ if $\eta(x)=0$ or $1$, and it is \emph{unstable} if $\eta(x)\geqslant 2$.
When a site is unstable it can \emph{topple}, sending away $2$ particles, each one to a neighbor of $x$ chosen independently at random.
Toppling a stable site is an \emph{illegal} operation and it will not be allowed.

To define the random topplings, start with an independent set of \emph{instructions} $\II=(\tau^{x,j}:x\in\Z^d, j\in\N)$, where $\tau^{x,j}=\tau_{xy}$ with probability $p(y-x)$.
$\PPP^\nu$ will denote the joint law of $\eta$ and $\II$, where $\eta$ has distribution $\nu$ and is independent of $\II$.

Let $\h=\big(h(x);x\in\Z^d\big)$ count the number of topplings at each site.
The toppling operation at $x$ is defined by $\Phi_x(\eta,\h)= \big(\tau^{x,2h(x)+2}\cdot\tau^{x,2h(x)+1}\cdot\eta, \h + \delta_x \big)$ and $\Phi_x\eta$ is a short for $\Phi_x(\eta,0)$.

The analogous representation for the ARW is described using the same notations.
For this model, a site $x$ is \emph{stable} in the configuration $\eta$ when $\eta(x)=0$ or $\eta(x)=\varrho$ and it is \emph{unstable} when $\eta(x)\geqslant 1$.
The instructions $(\tau^{x,j}:x\in\Z^d, j\in\N)$ are independent and are equal to $\tau_{xy}$ with probability $\frac{p(y-x)}{1+\lambda}$ or $\tau_{x\varrho}$ with probability $\frac{\lambda}{1+\lambda}$.
The ``toppling'' $\Phi_x(\eta,\h)= \big(\tau^{x,h(x)+1}\eta, \h + \delta_x \big)$ is \emph{legal} when $x$ is \emph{unstable} in $\eta$, i.e., when $\eta(x)\geqslant1$.

\paragraph{Properties.}
For $\alpha=(x_1,\dots,x_k)$, we write $\Phi_\alpha = \Phi_{x_k}\Phi_{x_{k-1}}\cdots\Phi_{x_1}$ and say that $\Phi_\alpha$ is \emph{legal} for $\eta$ if $\Phi_{x_l}$ is legal for $\Phi_{(x_{l-1},\dots,x_1)}\eta$ for each $1\leqslant l\leqslant k$.
In this case we say that $\alpha$ is a \emph{legal sequence}
of topplings for $\eta$.
Let $m_\alpha=\big(m_\alpha(x);x\in\Z^d\big)$ be given by
$m_\alpha(x)=\sum_l{\I}_{x_l=x}$, the number of times the site $x$ appears in
$\alpha$.
We write $m_\alpha \geqslant m_\beta$ if $m_\alpha(x) \geqslant m_\beta(x)\ \forall\ x$, and $\tilde{\eta} \geqslant \eta$ if $\tilde{\eta}(x) \geqslant \eta(x)\ \forall\ x$.
We also write $(\tilde{\eta},\tilde h) \geqslant (\eta,h)$ if $\tilde{\eta} \geqslant \eta$ and $\tilde{h} = h$.

Let $x$ be a site in $\Z^d$ and $\eta,\eta'$ be configurations.
\begin{property}
\label{property1dependenceonodometer}
If $\alpha$ is a legal sequence for $\eta$, then $\Phi_\alpha\eta$ depends on $\alpha$ only through $m_\alpha$.
\end{property}
\begin{property}
\label{property2localmonotone}
$\Phi_\alpha\eta(x)$ is non-increasing in $m_\alpha(x)$ and non-decreasing in $m_\alpha(z),z\ne x$.
\end{property}
\begin{property}
\label{property3stablemonotone}
If $x$ is unstable in $\eta$ and $\eta'(x)\geqslant\eta(x)$, then $x$ is unstable in $\eta'$.
\end{property}
\begin{property}
\label{property4preservedmonotone}
If moreover $\eta'\geqslant\eta$ then $\Phi_x\eta' \geqslant \Phi_x\eta$.
\end{property}

\begin{proof}
Property~\ref{property1dependenceonodometer}.
For sandpile-like models this follows simply from the abelianess of arithmetics, but the ARW has a slight subtlety.
For each site $y$ consider the action $\eta(y)\mapsto\tau\eta(y)$ of an operator $\tau$ on $\eta(y)$.
It is of the form $\tau_k:n\mapsto n+k$ for some $k\geqslant 0$ when $\tau=\tau_{xz}$, $x\ne y$, or of the form $\tau_{-1}:n\mapsto n-1$ when $\tau=\tau_{yx}$, or $\tau_{\varrho}:n\mapsto \varrho\cdot n$ when $\tau=\tau_{y\varrho}$, and the last two are legal operations only when $n\geqslant 1$.
Now if $\tau_{\varrho}$ is legal then
$\tau_{\varrho}\tau_{k}=\tau_{k}\tau_{\varrho}$ are legal, idem for $\tau_{-1}$,
and moreover $\tau_{k}\tau_{k'}=\tau_{k'}\tau_{k}$.
The $\tau_{-1}$ and $\tau_{\varrho}$ do not commute, but this does not pose any problem since the order of their appearance in $\Phi_{\alpha}$ is determined by the stack of instructions $\big(\tau^{y,j}\big)_{j\in\N}$ at $y$.

Property~\ref{property2localmonotone}.
For each site $x$, toppling sites other than $x$ cannot decrease $\eta(x)$, and
toppling site $x$ cannot increase $\eta(x)$.

Properties~\ref{property3stablemonotone}~and~\ref{property4preservedmonotone} are obvious.
\end{proof}

\paragraph{Consequences.}
Let $V$ be a finite subset of $\Z^d$.
A configuration $\eta$ is said to be \emph{stable} in $V$ if all the sites $x\in V$ are stable.
We say that $\alpha$ is \emph{contained} in $V$ if all its elements are in $V$, and we say that $\alpha$ \emph{stabilizes} $\eta$ in $V$ if every $x\in V$ is stable in $\Phi_\alpha\eta$.
\begin{lemma}[Least Action Principle]
\label{lemma1leastactionprinciple}
If $\alpha$ and $\beta$ are legal sequences of topplings for $\eta$ such that $\beta$ is contained in $V$ and $\alpha$ stabilizes $\eta$ in $V$, then $m_\beta \leqslant m_\alpha$.
\end{lemma}
\begin{remark*}
The sequence $\alpha$ need not be legal, it suffices that the topplings can be defined in some sense, and satisfy Properties~\ref{property1dependenceonodometer}--\ref{property4preservedmonotone}.
For an application that makes use of the Least Action Principle in such context, see~\cite{fey-levine-peres-10}.
\end{remark*}
\begin{proof}
We follow Eriksson~\cite{eriksson-96}.
The proof consists in observing that all topplings performed by $\beta$ are
necessary for stabilization.
Let $\beta$ be a legal sequence contained in $V$ and $m_\alpha \ngeqslant m_\beta$.
Write $\beta=(x_1,\dots,x_k)$ and $\beta^{(j)}=(x_1,\dots,x_j)$ for $j\leqslant k$.
Let $\ell=\max\{ j : m_{\beta^{(j)}} \leqslant m_{\alpha} \}<k$ and $y=x_{\ell+1}\in V$.
Now $y$ is unstable in $\Phi_{\beta^{(\ell)}}\eta$ since $\beta$ is legal, moreover $m_{\beta^{(\ell)}}\leqslant m_{\alpha}$ and $m_{\beta^{(\ell)}}(y)= m_{\alpha}(y)$ by definition of $\ell$.
By Properties~\ref{property2localmonotone}~and~\ref{property3stablemonotone}, $y$ is unstable for $\Phi_{\alpha}\eta$ and therefore $\alpha$ does not stabilize $\eta$ in $V$.
\end{proof}

\begin{lemma}[Abelian Property]
\label{lemma2abelianproperty}
If $\alpha$ and $\beta$ are both legal toppling sequences for $\eta$ that are contained in $V$ and stabilize $\eta$ in $V$, then $m_\alpha=m_\beta$.
In particular, $\Phi_\alpha\eta=\Phi_\beta\eta$.
\end{lemma}
\begin{proof}
Apply Lemma~\ref{lemma1leastactionprinciple} in two directions: $m_\alpha\leqslant m_\beta\leqslant m_\alpha$.
\end{proof}

By Lemma~\ref{lemma2abelianproperty}, $m_{V,\eta} = m_\alpha$ and $\xi_{V,\eta} = \Phi_\alpha\eta$ are well defined.

\begin{lemma}[Monotonicity]
\label{lemma3generalmonotonicity}
If $V\subseteq V'$ and $\eta \leqslant \eta'$, then $m_{V,\eta}\leqslant m_{V',\eta'}$.
\end{lemma}
\begin{proof}
Let $\alpha$ and $\beta$ be legal and stabilizing sequences respectively for $\eta'$ in $V'$ and for $\eta$ in $V$.
By successively applying Properties~\ref{property3stablemonotone}~and~\ref{property4preservedmonotone}, we see that
$\beta$ is also a legal sequence for $\eta'$ in $V$.
Thus $m_\alpha\geqslant m_\beta$ by Lemma~\ref{lemma1leastactionprinciple}
and the result follows.
\end{proof}

By monotonicity, the limit $m_\eta = \lim_{n}m_{V_n,\eta}$ exists and does not depend on the particular sequence $V_n\uparrow\Z^d$.
A configuration $\eta$ is said to be \emph{stabilizable} if $m_\eta(x)<\infty$ for every $x\in\Z^d$.\footnote
{If $\eta$ is stabilizable, $m_{V_n,\eta}(x) = m_\eta(x)$ for large $n$, thus by
Property~\ref{property2localmonotone} we have $\xi_{V_{n'+1}}(x) \geqslant
\xi_{V_{n'}}(x)$ for $n'>n$, and since $\xi_V(x)\leqslant 2$ for any $V$,
the final state $\xi(x)=\lim_{n\to\infty}\xi_{V_n}(x)$ is well defined and does
not depend on the particular sequence $V_n\uparrow\Z^d$.}

\begin{lemma}
\label{lemma4fixationstabilizable}
Let $\nu$ be a translation-invariant, ergodic distribution with finite density $\nu\big(\eta(0)\big)$.
Then $\Pb^\nu(\mbox{the system locally fixates})=\PPP^\nu\big(m_\eta(0)<\infty \big) \in \{0,1\}$.
\end{lemma}
\begin{proof}[Sketch of the proof]
The same proof works for both the ARW and the SSM, taking $\lambda=0$ for the latter.
We present the main steps and make references to~\cite{rolla-sidoravicius-11} when omitting details.
Let $h_t(x)$ denote the number of topplings at site $x$ during the time interval $[0,t]$, meaning any action performed at $x$, including unsuccessful attempts to sleep.
Write $h_\infty(x)=\lim_{t\to\infty}h_t(x)$, the limit exists as $h_t(x)$ is non-decreasing in $t$.

\medskip

The proof is split in three parts.
First, the $0$-$1$ law for the $\PPP^\nu$-probability that $\eta$ is stabilizable at some fixed site~$x$.
Second, that $\PPP^\nu$ provides a coupling for $\Pb^\nu$ so that
\begin{equation}
 \label{eq2fixationstabilizable}
 \Pb^\nu\big[h_\infty(x) \geqslant r\big]
 =
 \PPP^\nu\big[m_\eta(x)\geqslant r\big]
 \qquad
 \mbox{for each} \quad r>0.
\end{equation}
Third, that blowups do not happen in finite time:
\begin{equation}
 \label{eq3noblowups}
 \Pb^\nu\big[h_t(x) \geqslant r\big]\to 0 \quad\mbox{as}\quad r\to\infty
 \qquad\mbox{for each fixed }t.
\end{equation}
Let us show that these imply the lemma.
If $\PPP^\nu\big[m_\eta(x)<\infty\big]=1$, then it follows from~\eqref{eq2fixationstabilizable} that $\Pb^\nu\big[h_t(x)\mbox{ eventually constant}\big]=1$, thus $x$ is eventually stable in $\eta_t$ and in particular $\eta_t(x)$ remains bounded for large $t$.
But $\eta_t(x)$ can only decrease when $x$ is unstable, so $\Pb^\nu\big[\eta_t(x)\mbox{ converges}\big]=1$.
Otherwise $\PPP^\nu\big[ m_\eta(x)=\infty\ \big]=1$ by the $0$-$1$ law,
then~\eqref{eq2fixationstabilizable} gives $\Pb^\nu\big[h_t(x)\to\infty\mbox{ as
}t\to \infty\big]=1$ and by~\eqref{eq3noblowups} we know that
$\big(h_t(x)\big)_{t\geqslant 0}$ cannot blow up in finite time, whence for each
$x$, the value of $\eta_t(x)$ jumps for arbitrarily large times, and the system
stays active.

\medskip

The $0$-$1$ law comes from the fact that $m_{\eta}(x)=\infty$ implies
$m_{\eta}(x+z)=\infty$ for all $z$ such that $p(z)>0$, for a.e.\ $\II$.
See~\cite{rolla-sidoravicius-11} for the details, as well as the proof
of~\eqref{eq3noblowups}.

\medskip

The proof of~\eqref{eq2fixationstabilizable} is again divided in parts.
First we use $\PPP^\nu$ to simultaneously produce processes $(\eta_t^M,h_t^M)_{t\geqslant0}$
distributed as $\Pb^\nu_M$ for all $M\in\N$, starting from $\eta^M(z) = \eta(z) \I_{\|z\|\leqslant M}$ at $t=0$.
It follows from~\eqref{eq1finitebox} that
\begin{equation}
 \label{eq4limitmthent}
 \PPP^\nu\big(h_t^M(x)\geqslant r\big)
 =
 \Pb^\nu_M\big(h_t(x)\geqslant r\big)
 \mathop{\longrightarrow}_{M\to\infty}
 \Pb^\nu\big(h_t(x)\geqslant r\big)
 \mathop{\longrightarrow}_{t\to\infty}
 \Pb^\nu\big(h_\infty(x)\geqslant r\big).
\end{equation}
We check that
\begin{equation}
 \label{eq6limittthenm}
 \PPP^\nu\big(h_t^M(x)\geqslant r\big)
 \mathop{\longrightarrow}_{t\to\infty}
 \PPP^\nu\big(m_{\eta^M}(x)\geqslant r\big)
 \mathop{\longrightarrow}_{M\to\infty}
 \PPP^\nu\big(m_\eta(x)\geqslant r\big)
\end{equation}
and finally show that the
limits in~\eqref{eq4limitmthent} and~\eqref{eq6limittthenm} commute.

The coupling $\big\{(\eta_t^M,h_t^M)_{t\geqslant0}\big\}_{M\in\N}$ satisfying~\eqref{eq6limittthenm} is described in details in~\cite{rolla-sidoravicius-11}.
We show that the limits in~\eqref{eq4limitmthent} and~\eqref{eq6limittthenm} commute via monotonicity.
This follows from the coupling as well, for which $h_t^M(x)$ is a.s.\ non-decreasing in $M$ and $t$.
\end{proof}

\section{Strategy for stabilization}
\label{sec0heuristics}

In this section we give an overview of the general strategy in the proof of fixation.

The goal is to construct an algorithm that, given any pair $\eta$, $\II$, tries to stabilize all the particles initially present in $\eta$ following the jump instructions in $\II$, but with the aid of some extra topplings.
The algorithm must be such that $m_\eta(0)=0$ whenever it is successful.
It may well fail, but all we need is \emph{a strategy that succeeds with positive probability}, which in turn implies a.s.\ fixation by the $0$-$1$ law (see Lemma~\ref{lemma4fixationstabilizable}).

The algorithm consists of applying a \emph{settling procedure} to each particle.
This procedure \emph{explores} $\II$ until it identifies a suitable \emph{trap} for the particle.
The exploration follows the path that the particle would perform if we always toppled the site it occupies, and stops when the trap has been chosen.
In the absence of a suitable trap, we declare the algorithm to have failed.
The trap always lies on the exploration path, however not necessarily on its tip because we need to explore further away before taking the decision.
Once the trap has been chosen, the particle is moved along the exploration path
until it reaches the trap, where it is settled.

\medskip

The major issue is of course the spread of activity.
Thus the first requirement for the settling procedure is \emph{not to disturb particles that have already been settled}, meaning that the exploration does not examine those sites.

We need to conciliate this requirement with the goal that the algorithm have positive probability of success.
This imposes that the settling procedures not only succeed with a high probability, but these probabilities should actually converge to $1$.
We hence consider a procedure that tries to find the traps close together as much as possible, in order to \emph{leave more space for the maneuver of subsequent steps}.

Moreover, certain control is needed on the joint distribution of the outcome of different explorations.
Some of the explored instructions are actually not going to be used by the
corresponding particle by the time it settles at the trap, leaving some
\emph{corrupted sites} that may interfere with the next steps.
In our proof we go for \emph{independence}, which means that \emph{corrupted sites left by previous steps must be avoided}.
Hence the settling procedure should try not only to find its trap close to those already found, but also to \emph{keep the corrupted region as compact as possible}.

\medskip

We seek a procedure that sets the trap close to the positions where the previous
traps have been placed, therefore implying typically long excursions away from
the particle starting position.
But \emph{we are bound to follow the instructions in $\II$,} and cannot prevent the particles from stopping before.

To circumvent this restriction, we \emph{enforce activation and push the
particles further}, in a way that $\tilde m$, the total number of topplings
after enforcement, provides an \emph{upper bound} for the true $m_\eta$.
For the ARW this is done by ignoring some instructions of the form
$\tau_{x\varrho}$, and we show in the next session that it increases $m_\eta$.
The SSM case is considered in the sequel.

\medskip

Let us describe the stabilizing strategy used in this paper.
Label the initial positions of the particles on $\Z_+$ by
$0 \leqslant x_1 \leqslant x_2 \leqslant x_3 \leqslant \cdots$.
Take $a_0=0$ and suppose the first $k-1$ traps have been successfully set up at
positions $0<a_1<a_2<\cdots<a_{k-2}<a_{k-1}<x_{k-1}$.
The set $\llbracket 0,a_{k-1} \rrbracket$ contains all traps and corrupted sites
found so far,
where $\llbracket z,w \rrbracket$ denotes $[z,w]\cap\Z$.

The settling procedure starts with an \emph{exploration}.
Starting at $x_k$, examine and follow the instructions in $\II$ one by one
(whenever a sleep instruction is found, the next instruction
at the same site is to be examined).
Follow this exploration until reaching $a_{k-1}$.

Next we \emph{set up the trap}.
We consider first the ARW and postpone the SSM subtleties.
During the $k$-th exploration, we are sure to visit every $x\in
B_k=\llbracket a_{k-1}+1,x_k-1 \rrbracket$.
Moreover, the last instruction explored at each $x\in B_k$ is a jump to the
left, see Figures~\ref{fig1arw}~and~\ref{fig2arw}.
For some $x\in B_k$, the second last instruction may be a sleep instruction.
Let $a_k$ be the leftmost such site.
If there is none, we declare the procedure \emph{unsuccessful} and stop.

At this point the instruction manipulation enters the game: we suppress all the
sleep instructions explored meanwhile, except for the last sleep instruction at
$a_k$, which becomes the trap.

Since the trap is a sleep instruction found immediately before the last
instruction, which is a jump to the left, we are sure that the exploration path
has not gone to the right of $a_k$ after exploring the trap, so \emph{all the
corrupted sites will be in $\llbracket a_{k-1}+1,a_k\rrbracket$}, see
Figure~\ref{fig2arw}.
We \emph{repeat this procedure indefinitely}, unless it fails at some step.

We finally show that this strategy is \emph{successful with positive
probability}.
Each site $x \in B_1$ has a sleep instruction just before its last jump
independently and with probability $\bar\mu=\frac{\lambda}{1+\lambda}$.
Thus, $a_1-a_0$ is distributed as a geometric random variable of mean
$\bar\mu^{-1}$, truncated at $x_1-a_0$.
Since no corrupted sites were left outside $\llbracket a_0+1,a_1 \rrbracket$,
the interdistance $a_2-a_{1}$ is independent of $a_1$ and has the same
distribution.
By the law of large numbers $a_n \sim n/\bar\mu$.
On the other hand, $x_n \sim n/\mu$.
Therefore, if $\mu<\bar\mu$ there is positive probability that $a_k < x_k$ for
all $k$, which implies success.

\medskip

In the SSM case, an isolated particle stands still and the site where it is
located is stable.
A natural way to enforce activity is to allow a single particle to move.
To achieve this, we regard the toppling operation as having two
distinct steps: the particles are sent to neighboring sites, but one at a time,
leading to what we call \emph{half-topplings}.

However, this has to be consistent with the specification of the model.
More precisely, to determine whether a half-toppling is legal we need to know
whether or not the site has been half-toppled.
If $\eta(x) \geqslant 2$ the half-toppling is always legal.
It is also legal if $\eta(x) = 1$, but only when $x$ has been half-toppled.
These rules yield the same $m_\eta$, and moreover, violating the latter
condition, i.e., half-toppling stable particles, results in an upper
bound.
See Section~\ref{sec5phasetransitionssm} for the details of how this tool is
used in the proof of fixation.

\section{Phase transition for activated random walks}
\label{sec1phasetransition}

In this section we prove Theorem~\ref{theo1fixation}.
The core of the proof of fixation is an algorithm that, given $\II$ and
$\eta$, selects
a subset of the sleep instructions in $\II$, providing a new $\tilde \II$ under
which $\tilde m_{\eta}(0)=0$.
This algorithm works with positive probability, and the theorem follows from
Lemma~\ref{lemma4fixationstabilizable} and Lemma~\ref{lemma1wakeupmonotonicity}
below.

Besides the $\tau_{xy}$ and $\tau_{x\varrho}$ considered in
Section~\ref{sec1flag}, consider in addition the neutral instruction $\iota$,
given by $\iota\eta=\eta$.
Replacing a sleep instruction by a neutral instruction may be seen as enforcing
the active state on a passive particle, or replacing the transition $A\to S$ by
$A\to S\to A$.
Given two sets of instructions $\II=(\tau^{x,j})_{x,j}$ and
$\tilde\II=(\tilde\tau^{x,j})_{x,j}$, we write $\tilde \II\geqslant \II$ if for
every $x\in\Z^d$ and $j\in\N$, either $\tau^{x,j}=\tilde \tau^{x,j}$, or
$\tau^{x,j}=\tau_{x\varrho}$ and $\tilde \tau^{x,j}=\iota$.
As we deal with more than one set of instructions, we need to enlarge the
notations of Section~\ref{sec1flag} and write $m_{V,\eta;\II}$,
$\Phi_{\alpha;\II}$, etc., to specify which set of instructions $\II$ is being
considered.

\begin{lemma}[Monotonicity with enforced activation]
\label{lemma1wakeupmonotonicity}
Let $\tilde\II\geqslant\II$ be sets of instructions, $\eta$ be a initial state,
and $V\subseteq \Z^d$ be a fixed finite set.
Then $m_{V,\eta;\II} \leqslant m_{V,\eta;\tilde \II}$.
In particular, $m_{\eta;\II} \leqslant m_{\eta;\tilde \II}$.
\end{lemma}
\begin{proof}
Let $\tilde \II \geqslant \II$.
By Lemma~\ref{lemma1leastactionprinciple}, it suffices to show that if $(\tilde
\eta,\tilde h) \geqslant (\eta,h)$, and $\alpha$ is a legal sequence for
$(\eta,h)$
under $\II$, then $\alpha$ is also legal for $(\tilde \eta,\tilde h)$ under
$\tilde \II$.
We show the stronger fact that $\Phi_{\alpha;\tilde{\II}}(\tilde{\eta},\tilde
{h}) \geqslant \Phi_{\alpha;{\II}}({\eta},{h})$.
The proof is by induction on the number of steps $k=\sum_z m_\alpha(z)$.

For $k=0$ the claim is trivially true.
For $k \geqslant 2$ we write $\alpha=\alpha'\alpha''$ and apply the induction
hypothesis for both $\alpha'$ and $\alpha''$.
It remains to consider $k=1$.
Write $\alpha = (y)$.
Since $\tilde\eta\geqslant\eta$, $y$ is also unstable in $\tilde\eta$.
Write $\eta'=\tau^{y,h(y)+1}\eta$ and $\tilde\eta' =
\tilde\tau^{y,h(y)+1}\tilde\eta$.
Now either $\tau^{y,h(y)+1}=\tilde\tau^{y,h(y)+1}$, and thus $\tilde
\eta'=\tau^{y,h(y)+1}\tilde \eta \geqslant \tau^{y,h(y)+1}\eta=\eta'$, or
$\tau^{y,h(y)+1}=\tau_{x\varrho}$ and $\tilde \tau^{y,h(y)+1}=\iota$, and thus
$\eta'=\tau^{y,h(y)+1}\eta = \tau_{y\varrho}\eta\leqslant \eta\leqslant
\tilde\eta = \iota \tilde\eta=\tilde\eta'$.
In any case, $\eta'\leqslant\tilde\eta'$ and
$\Phi_{y;\tilde{\II}}(\tilde{\eta},\tilde {h})=(\tilde\eta',h+\delta_y)
\geqslant (\eta',h+\delta_y)=\Phi_{y;{\II}}({\eta},{h})$.
\end{proof}

\begin{proof}[Proof of Theorem~\ref{theo1fixation}]
Let $\PPP^\mu$ denote $\PPP^\nu$ for $\nu$ the product probability distribution having Poisson($\mu$) as marginals.
By Lemmas~\ref{lemma3generalmonotonicity} and~\ref{lemma4fixationstabilizable}
it suffices to show that $\PPP^\mu\big(m_\eta(0)=0\big)>0$ when $\mu<\frac
\lambda{1+\lambda}$ and $\PPP^\mu\big(m_\eta(0)=\infty\big)>0$ when $\mu>1$.

\medskip

Let $\mu>1$.
For some $\delta>0$, the $\PPP^\mu$-probability that there are at least $\mu M$ particles in $V=[-M,0]$ is at least $2\delta$, independently of $M$.
After stabilizing $\eta$ in $V$, at least $(\mu-1)M$ particles will have to
exit $V$, and because of the topological constraints of the one-dimensional
lattice, $m_{V,\eta}(-M)\geqslant (\mu-1)M/2$ or
$m_{V,\eta}(0)\geqslant (\mu-1)M/2$ must happen.
Thus, at least one of these two events has probability bigger than $\delta$.
The latter case gives directly $\PPP^\mu\big(m_\eta(0)\geqslant
(\mu-1)M/2\big)>\delta$, and in the former case we can consider $V=[0,M]$, which
gives the same inequality by translation invariance.
Letting $M\to\infty$ we get
$\PPP^\mu\big(m_\eta(0)=\infty\big)\geqslant\delta>0$.

\medskip

In the remainder of this section we prove that $\eta$ is stabilizable
when $\mu<\frac \lambda{1+\lambda}$.
Given $\II$ and $\eta$, the goal is to construct a stabilizing strategy in an
algorithmic way, finding $\tilde\II\geqslant\II$ such that
$m_{\eta;\tilde\II}(0)=0$.
We consider first the recurrent case $p(-1)=p(+1)=\frac12$.

\medskip

Start placing an imaginary \emph{barrier} at site $a_0=0$ and pick the first
particle to the right of this barrier.
\emph{Explore} the instructions in $\II$ one by one, following a path
that starts at the position $x_1$ of this first particle.
This means following the path that would have been performed by this particle if we were to topple the site that contains it over and over.
Follow this exploration until it reaches the origin.
When sleep instructions are found, this explorer looks for the next instruction
at the same site,
until it finds a jump instruction, as in Figure~\ref{fig1arw}.

\medskip

We then \emph{set up the trap} for the first particle.
For each site $y\in \llbracket a_0+1,x_1-1 \rrbracket$, the last instruction
examined by this
explorer at $y$ must be a jump to the left, $\tau^{y,j}=\tau_{yz}$,
where $z=y-1$ and $j=j_1(y)$ is the total number of instructions examined at
site $y$.
It may be the case that before this jump to the left, the previous instruction $\tau^{y,j-1}$ found at $y$ is $\tau_{y\varrho}$.
Let
$a_1=\min\big\{y\in \llbracket a_0+1,x_1-1 \rrbracket
:\tau^{y,j_1(y)-1}=\tau_{y\varrho}\big\}$ be
the leftmost such site.
This is the site where the \emph{trap is set up}.
Place another imaginary \emph{barrier} at $a_1$, as in Figure~\ref{fig1arw}.
\begin{figure}[!htb]
 \centering
 \includegraphics[scale=.5]{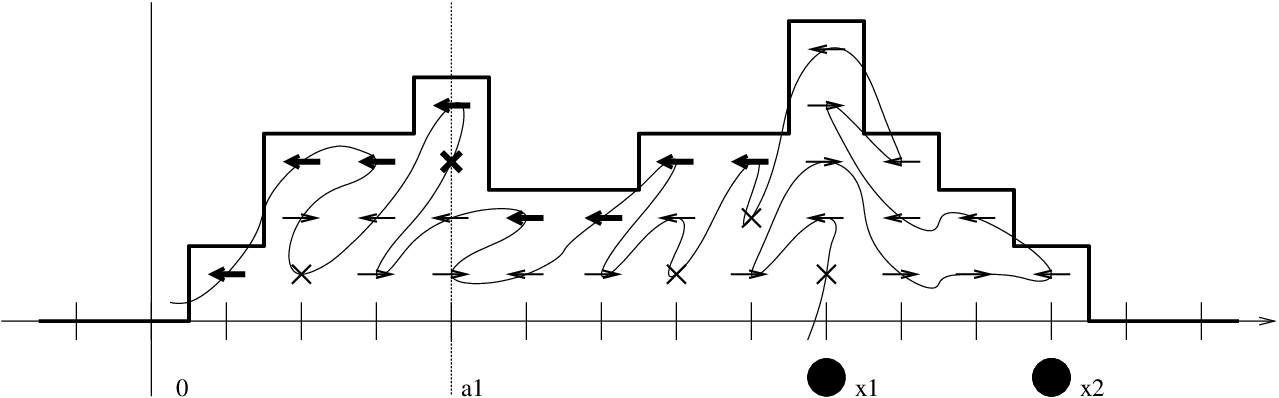}
 \caption{
 First exploration path for the ARW.
 It starts at position $x_1$ of the first particle and stops when it reaches the origin.
 The horizontal axis represents the lattice, and above each site $x$ there is a sequence of instructions $(\tau^{x,j})_j$.
 The bold arrows indicate the last jump found at each site $x \in \llbracket
1,x_1-1 \rrbracket$, and the bold cross indicates a sleep instruction found just
before the last jump, this being the leftmost such cross, whose
location defines $a_1$.
 }
 \label{fig1arw}
\end{figure}
If no such site is found, we declare the procedure to have failed and we stop
it.
Otherwise, we declare the first step to be \emph{successful}.

\medskip

We now \emph{explore the path} for the \emph{second particle} and \emph{set up
the next trap}.
Start another explorer at $x_2$, the position of the second particle to the
right of the origin, and follow this exploration until hitting the barrier
$a_1$.
This explorer will examine the instructions of $\II$ in the order as they
appear, but skipping those that have been used in the previous step, as in
Figure~\ref{fig2arw}.
\begin{figure}[!htb]
 \centering
 \includegraphics[scale=.5]{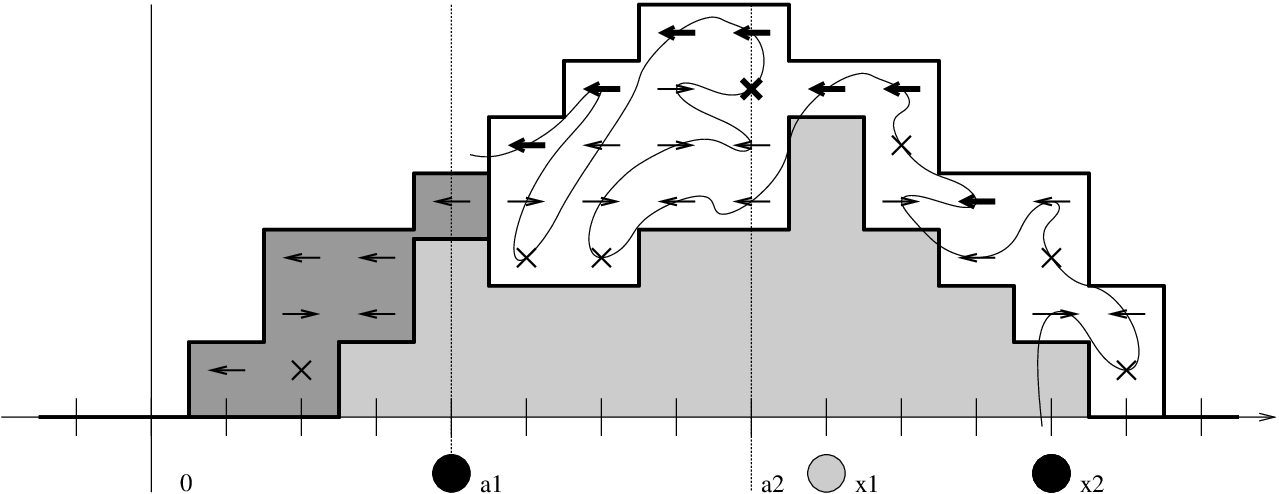}
 \caption{
 Second exploration path for the ARW.
 It starts at position $x_2$ of the second particle and stops when it reaches $a_1$.
 The regions in gray indicate the instructions already examined by the first
explorer.
 The dark gray contains instructions examined but not used, whose locations
determine the set of corrupted sites.
 }
 \label{fig2arw}
\end{figure}
As before, for each $y\in \llbracket a_1+1,x_2-1 \rrbracket $, the last
instruction at $y$
must be a jump to the left, $\tau^{y,j}=\tau_{yz}$, where $z=y-1$ and $j=j_2(y)$
is the total number of instructions examined at site $y$ by both explorers.
The \emph{trap} is set up at $a_2=\min\big\{y\in\llbracket a_1+1,x_2-1
\rrbracket: \tau^{y,j_2(y)-1}=\tau_{y\varrho}\big\}$, the leftmost site where
the
instruction examined just before the last jump to the left is a sleep
instruction, as in Figure~\ref{fig2arw}.
Again we place a new imaginary \emph{barrier} at $a_2$.
If no such site is found, we declare the procedure to have failed and we stop
it.
Otherwise, we declare the second step to be \emph{successful}.

\medskip

This procedure can be carried on indefinitely, as long as \emph{all the steps
are successful}.
We then perform a similar construction for the negative half line, finding the
site $x_{-1}$ with the first particle to the left of the origin, setting up a
trap at $a_{-1}\in \llbracket x_{-1}+1,-1 \rrbracket $, then at
$a_{-2}\in \llbracket  x_{-2}+1,a_{-1}-1 \rrbracket $, and so on.

\medskip

Let us show that \emph{success implies stability}.
Suppose all the steps in the algorithm are successful.
Replace each sleep instruction in $\II$ by a neutral instruction, \emph{except
for those providing the traps}.
Namely, take
\[
  \tilde\tau^{y,j} =
  \begin{cases}
    \iota, & \tau^{y,j}=\tau_{y\varrho},
            (y,j)\ne\big(a_{n},j_n(a_n)\big)\ \forall n\in\Z, \\
    \tau^{y,j}, &\mbox{otherwise.}
                     \end{cases}
\]
Fix $n\in\N$.
We need to show that, following the instructions of $\tilde\II$, $\eta$ is
stabilized in $V_n=[x_{-n},x_n]$ with finitely many topplings, and moreover
$m_{\eta;\tilde\II}(0)=0$.

We first stabilize the particle that starts at $x_1$.
To that end, we successively topple the sites found by the first explorer.
Since all the sleep instructions have been replaced, the trajectory of the
particle started at $x_1$ will coincide with the path of the corresponding
explorer, until it finds the trap
$\tau^{a_1,j_1(a_1)-1}=\tau_{a_1\varrho}$.
At this moment the particle will become passive, and the site $a_1$ will be stable.

Notice that, after the last visit to $a_1$, the explorer does not go
further to the right, so when settling the first particle we use all the
instructions examined so far, except some lying in $\llbracket a_0+1,a_1
\rrbracket$.
Therefore, the same procedure can be applied to the second particle, as it will
find the same instructions that determined the second exploration path.

Notice also that the first particle does not visit $0$, and the second particle
neither visits $0$ nor $a_1$, thus it is settled without activating the first
particle.
Following the same procedure, the $k$-th particle is settled at $a_k$, without ever visiting $\{0,a_1,a_2,\dots,a_{k-1}\}$, for all $k=1,\dots,n$.
After settling the $n$ first particles in $\Z_+$, we perform the analogous procedure for the first $n$ particles in $\Z_-$.

This means that $\eta$ can be stabilized in $V_n$ with finitely many topplings, not necessarily in $V_n$, and never toppling the origin.
By Lemma~\ref{lemma1leastactionprinciple}, $m_{V_n,\eta;\tilde \II}(0)=0$.
On the other hand, it follows from Lemma~\ref{lemma1wakeupmonotonicity} that
$0\leqslant m_{V_n,\eta;\II}(0)\leqslant m_{V_n,\eta;\tilde \II}(0)$, whence
$m_{V_n,\eta;\II}(0)=0$.
Since it holds for all $n\in\N$ and $V_n\uparrow\Z$ as $n\to\infty$, this gives
$m_{\eta;\II}(0)=0$.

\medskip

We conclude with the proof that there is \emph{positive probability of success}
for the algorithm described above.
Namely, we show that the set of $(\eta,\II)$ for which the above construction is
successful has positive $\PPP^\mu$-probability whenever
$\mu<\frac\lambda{1+\lambda}$.

We claim that the position $a_1$ of the first barrier has a geometric
distribution with parameter $\frac\lambda{1+\lambda}$.
To be more accurate, the claim is that the probability space can be enlarged so that we can define a random variable $Y_1$, independent of $\eta$, satisfying $\PPP^\mu\big(Y_1>k\big)=\big(\frac1{1+\lambda}\big)^k$, with the property that the first step of the construction is successful if and only if $Y_1 < x_1$, in which case the position $a_1$ of the first barrier is given by $a_1=Y_1$.

Let us prove the above claim.
For each site $y\in\Z$, the instruction $\tau^{y,j}$ can be a sleep instruction $\tau_{y\varrho}$ with probability $\frac{\lambda}{1+\lambda}$ or a jump instruction of the form $\tau_{yz}$ with probability $\frac{1}{1+\lambda}$, independently for each $j\in\N$.
Conditioning on $J^y=\{j^1<j^2<\cdots\}$, the subset of $\N$ for which $\tau^{y,j}$ is a jump instruction, $\tau^{y,j^k}=\tau_{yz}$ with probability $p(z-y)$, independently for each $k$.
Now the trajectory performed by the first explorer depends only on the
jump instructions, and, conditioned on the trajectory, we have
$\tau^{y,j_1(y)-1}=\tau_{y\varrho}$ with probability $\frac{\lambda}{1+\lambda}$
and independently for each $y\in \llbracket a_0+1,x_1-1 \rrbracket$.
If it fails for all such $y$, which happens with probability $\big(\frac1{1+\lambda}\big)^{x_1-1}$, sample $Y_1$ as $\PPP^\mu\big(Y_1>k\big|Y_1>x_1-1\big)=\big(\frac1{1+\lambda}\big)^{k-x_1+1}$.
Otherwise take $Y_1$ as the smallest such $y$.
This implies that $Y_1$ has the prescribed distribution, proving the claim.

Since each explorer skips the instructions already examined at previous steps,
the sequence $(a_1, a_2-a_1,a_3-a_2, \dots)$ is i.i.d.
More precisely, by the same argument as above, there is a sequence of i.i.d.\
variables $Y_1,Y_2,Y_3,\dots$ with the property that the $n$-th step is
successful if and only if the previous steps are successful and
$a_{k-1}+Y_k<x_k$, in which case $a_k=a_{k-1}+Y_k$.

Now the expectation of $Y_1$ is $\frac{\lambda+1}{\lambda}$.
By the law of large numbers, the putative position $X_k=Y_1+\dots+Y_k$ of the $k$-th barrier grows like $\frac{\lambda+1}{\lambda}k$.
On the other hand, by the law of large numbers, the initial position $x_k$ of the $k$-th particle to the right of the origin grows like $k/\mu$.
Therefore, as $\mu<\frac\lambda{1+\lambda}$, $X_k < x_k$ for all $k$ with
positive probability.

When we apply the same construction switching left and right,
we define $X_k=Y_k+\dots+Y_{-1}$ which satisfy $x_k < X_k$ for all
$k<0$ with the same probability and independently of the first part.
Therefore the construction is successful with positive probability and the proof of
the recurrent case is finished.

\medskip

We now consider the transient case $p(+1)=1-p(-1)>\frac12$.
If the $n$-th explorer never hits $a_{n-1}$, it is a.s.\ the case
that it examines only finitely many instructions at each site.
In this case we take $a_n=a_{n-1}$, no sleep instruction is needed for the
$n$-th particle as its trap will be at $\infty$, and it is still possible to
define the subsequent explorers.
When stabilizing the configuration $\eta$ in $V_n=\llbracket x_{-n},x_n
\rrbracket$, we let the
particles follow the path of the corresponding explorer as before.
The total number of topplings may be infinite (and in the limit the transient particles disappear), but each site is toppled finitely many times and, as in the symmetric case, the origin never topples.
\end{proof}

\section{Phase transition for the stochastic sandpile model}
\label{sec5phasetransitionssm}

In this section we prove Theorem~\ref{theo1fixationssm}.
The proof of fixation consists in an algorithm that, given $\II$ and $\eta$,
stabilizes $\eta$ with a sequence of topplings, \emph{some of which are not
legal}.
More precisely, we introduce the half-toppling operation, which moves
the particles one by one.
The configuration is then stabilized with the aid of some illegal
half-topplings, and without toppling the origin.
This algorithm works with positive probability, and the theorem follows from
Lemma~\ref{lemma4fixationstabilizable} and
Lemma~\ref{lemma1halftopplingmonotonicity} below.

A \emph{half-toppling} consists of sending
only one particle to a neighboring site chosen at random (i.e., following the
appropriate instruction) and is possible when $\eta(x)>0$.
We say that site $x$ is \emph{stable} if
$\eta(x)=0$ and \emph{unstable} if $\eta(x)\geqslant 2$.
When $\eta(x)=1$, the site $x$ can be stable or unstable, depending on whether
it has been half-toppled an even or an odd number of times.
A half-toppling at $x$ is considered \emph{legal} when $x$ is unstable, and is
considered \emph{semi-legal} when $\eta(x)\geqslant 1$ regardless of the site
being stable or not.
We denote a half-toppling at $x$ by $\phi_x$, that is, $\phi_x(\eta,\h)=
\big(\tau^{x,2h(x)+1}\eta, \h + \frac12\delta_x \big)$ and
$\phi_x\eta=\phi_x(\eta,0)$.
Notice that $\Phi_x = \phi_x\circ\phi_x$.

For a sequence $\beta=(y_1,\dots,y_n)$, define $\tilde m_\beta(y)=
\sum_j{\I}_{y_j=y}$ as the number of times that the site $y$ appears in
$\beta$.
Properties~\ref{property1dependenceonodometer}--\ref{property4preservedmonotone}
also hold with half-topplings.
For $\alpha=(x_1,\dots,x_k)$, denoting
$\alpha^2=(x_1,x_1,x_2,x_2,\dots,x_k,x_k)$, we have that $\tilde{m}_{\alpha^2}
= 2m_\alpha$ and, given a configuration $\eta$, the sequence of half-topplings
$\phi_{\alpha^2}$ is legal if and only if the sequence of topplings
$\Phi_\alpha$ is itself legal, and in this case
$\phi_{\alpha^2}\eta=\Phi_\alpha\eta$.

\begin{lemma}[Least Action Principle, with half-topplings]
\label{lemma1halftopplingmonotonicity}
Let $\beta$ be a sequence of half-topplings in the volume $V$ that is legal for
the configuration $\eta$.
Let $\alpha$ be a sequence of half-topplings in $\Z^d$ that is semi-legal for
$\eta$ and stabilizes $\eta$ in $V$.
Then $\tilde m_\beta \leqslant \tilde m_\alpha$.
\end{lemma}
\begin{proof}
The proof is the same as for Lemma~\ref{lemma1leastactionprinciple}.
\end{proof}
\begin{cor}[Abelian Property, with half-topplings]
\label{cor1abelianwithhalftopplings}
If $\alpha$ is a sequence of half-topplings in $V$ that is legal for $\eta$ and
stabilizes $\eta$ in $V$, then $\tilde m_\alpha= 2 m_{V,\eta}$.
\end{cor}
\begin{proof}
Take a sequence of topplings $\beta$ in $V$ that is legal for $\eta$ and stabilizes $\eta$ in $V$.
Then $\beta^2$ is a legal sequence of half-topplings in $V$ that is legal for
$\eta$ and stabilizes $\eta$ in $V$, and $\tilde m_{\beta^2} = 2 m_\beta$.
Now $m_\beta=m_{V,\eta}$ by Lemma~\ref{lemma2abelianproperty} and $\tilde m_{\alpha} \leqslant \tilde m_{\beta^2} \leqslant \tilde m_{\alpha}$ by Lemma~\ref{lemma1halftopplingmonotonicity}.
\end{proof}

\begin{proof}[Proof of Theorem~\ref{theo1fixationssm}]
The outline of the proof is the
following.
First we show that there is no fixation for $\mu>1$.
Next we describe how to stabilize the particles initially present in $\Z_+$,
and show that this scheme has positive probability of success if $\mu$ is small
enough.
We then argue by symmetry that an analogous procedure works for the particles
starting in $\Z_-$.
The first and last parts are identical to those presented in
Section~\ref{sec1phasetransition}, and we restrict this section to the
specificity of the settling procedure.

\medskip

Given $\II$ and $\eta$, we will construct a stabilizing strategy in an
algorithmic way, giving  $\big(\tilde m_{V_n}(x):x\in\Z\big)_{n=1,2,\dots}$,
where $\tilde m_{V_n}$ counts the number of semi-legal half-topplings (not
necessarily legal nor contained in $V_n$) that stabilize $\eta$ in $V_n$ and
satisfy $\tilde m_{V_n}(0)=0$, for some sequence $V_n \uparrow \Z_+$.

\medskip

Start placing an imaginary \emph{barrier} at site $a_0=0$ and pick the first
particle to the right of this barrier.
\emph{Explore} the instructions in $\II$ one by one, starting at the position
$x_1$ of this first particle and following its path that would have been
performed by this particle if we were to half-topple the site that contains it
over and over until it reached $a_0$, as in Figure~\ref{fig1grain}.
\begin{figure}[!htb]
 \centering
 \includegraphics[scale=.5]{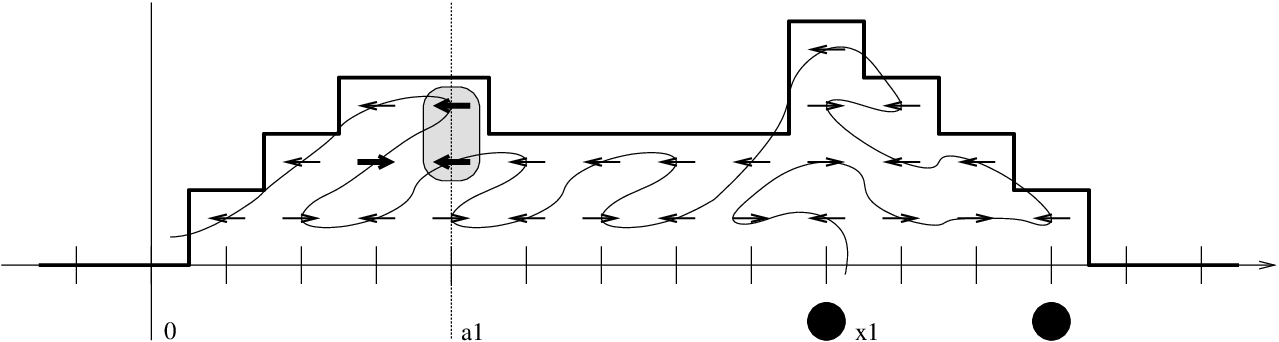}
 \caption{
 First exploration for the SSM.
 It starts at position $x_1$ of the first particle and stops when it reaches the origin.
 The horizontal axis represents the lattice, and above each site $x$ there is a sequence of instructions $(\tau^{x,j})_j$.
 The bold arrows represent the last jump to the right of this explorer and
the two last jumps (necessarily to the left) on its neighboring site $a_1$.
 }
 \label{fig1grain}
\end{figure}

We then \emph{set up the trap} for the first particle.
Fix the site
$a_1 \in \llbracket a_0+1,x_1 \rrbracket $ to the right of where the
explorer does its \emph{last right jump}.
The \emph{trap} will be set at $a_1$, as we explain later on.
Place a new barrier at this site, as
in Figure~\ref{fig1grain}.
If no such site is found, we declare the procedure to have failed and stop it.
Otherwise, we declare the first step to be \emph{successful}.

Each step of the explorer corresponds either to a semi-legal half-toppling of a
stable site, or a legal half-toppling of an unstable site.
If a site is visited more than once, one of the last two moves necessarily
corresponds to the former case.
Now notice that $a_1$ is visited at least twice by the respective explorer
(which implies that at least one of the last two visits corresponds to an illegal
toppling) and moreover, no site to the right of $a_1$ is visited after the two
last jumps.

\medskip

We now \emph{explore the path} for the \emph{second particle} and \emph{set up
the next trap}.
Start another exploration at $x_2$, the position of the second particle to the
right of the origin, and follow this exploration until it hits the barrier
$a_1$.
This explorer will examine the instructions of $\II$ in the order as they
appear, but skipping those that have been used in the previous step, as in
Figure~\ref{fig2grain}.
\begin{figure}[!hbt]
 \centering
 \includegraphics[scale=.5]{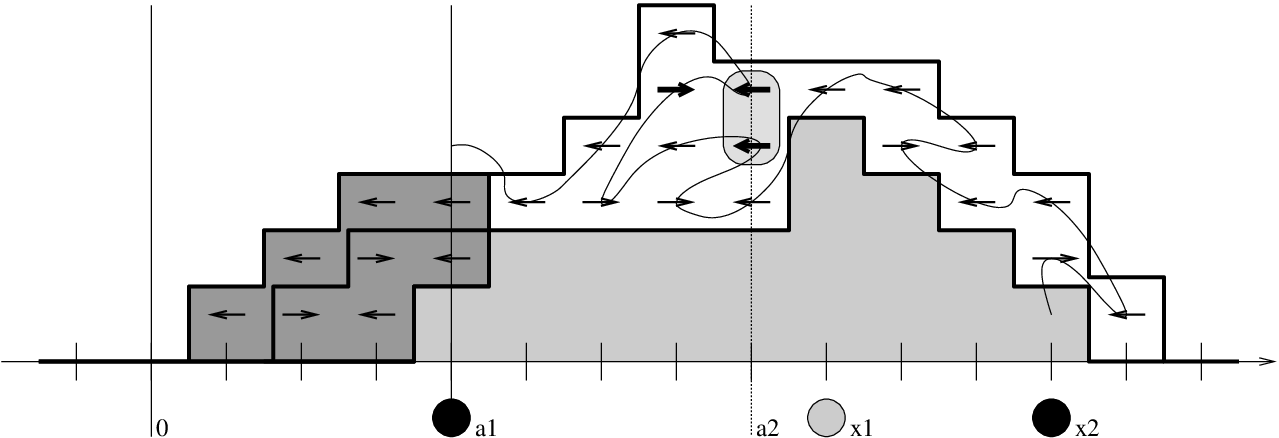}
 \caption{
 Second exploration for the SSM.
 It starts at position $x_2$ of the second particle and stops when it reaches $a_1$.
The regions in gray indicate the instructions already examined by the first
explorer.
The dark gray contains instructions examined but not used, whose locations
determine the set of corrupted sites.
A subset of the dark gray instructions may be used, depending on how the first
particle is settled at $a_1$.
}
\label{fig2grain}
\end{figure}
As before, fix the site $a_2\in \llbracket a_1+1, x_2\rrbracket $ to the right
of where this explorer does its last right jump.
The \emph{trap} will be set at $a_2$, and we place a new barrier at this site,
as in Figure~\ref{fig2grain}.
If no such site is found, we declare the procedure to have failed and stop it.
Otherwise, we declare the second step to be \emph{successful}.
This procedure can be carried on indefinitely, as long as \emph{all the steps
are successful}.

\medskip

We now show that \emph{success implies stability}.
Suppose all the steps in the algorithm are successful.
Let us show how to stabilize $\eta$ in $\llbracket 0,x_n \rrbracket $.
We first stabilize the particle that starts at $x_1$.
To that end, we successively half-topple the sites found by the first explorer,
except for the second last half-toppling at $a_1$.
At this point, if $\tilde m(a_1)\in 2\N_0$, the site $a_1$ is stable and we
\emph{abstain from further illegal half-topplings}, leaving the particle at
$a_1$.
Otherwise, if $\tilde m(a_1)\in 2\N_0+1$, the site is unstable, in which case we
keep the half-topplings until the last half-toppling at $a_1$, which must be
stable at this time, again leaving the particle at $a_1$.

Notice that, after each of the two last visits to $a_1$, the explorer does
not go further to the right.
Thus, when settling the first particle we use all the
instructions examined by the first explorer, except some lying in
$\llbracket a_0+1,a_1 \rrbracket$.
Therefore, this procedure can be applied to the second particle as it will find
the same instructions that determined the exploration path.

Notice also that the first particle does not visit $0$, and the second particle neither visits $0$ nor $a_1$, so it is settled without causing further instabilities.
Following the same procedure, the $k$-th particle is settled at $a_k$, without ever visiting $\{0,a_1,a_2,\dots,a_{k-1}\}$, for all $k=1,\dots,n$.

\medskip

To conclude the proof of the theorem, we need to show that the construction is
\emph{successful with positive probability} when $\mu<\frac14$.

We claim that the position $a_1$ of the first barrier is distributed as the time of the first jump to the left of a discrete-time simple symmetric random walk $(S_n)_{n=0,1,2,\dots}$ started at $S_0=0$ and conditioned to stay positive forever after.
To be more accurate, the claim is that the probability space can be enlarged so
that we can define a random variable $Y_1$, independent of $\eta$, distributed
as $\inf\{n:S_{n+1}=S_n-1\}$, with the property that the first step of the
construction is successful if and only if $Y_1 < x_1$, in which case the
position $a_1$ of the first barrier is given by $a_1=Y_1$.

Let us prove the above claim.
A visual proof is contained in Figure~\ref{fig1grain}.
Let $(\tilde S_n)_{n=0,1,\dots,k}$ denote the path of the first explorer.
If $(S_n)_{n=0,1,\dots,k}$ denotes the reversed path $S_n = \tilde S_{k-n}$, then $S_n$ satisfies $S_0=0$, $S_k=x_1$ and $S_n>0$ for all $n>0$.
Notice that each $(S_n)_n$ with the above properties is possible, and it happens with probability proportional to $2^{-k}$, that is, $2^{-k}/Z$, where $Z$ is the sum over all $k'$ of $2^{-k'}$ times the number of paths $(S_n)_{n=0,1,\dots,k'}$ with the above properties.
But this is the probability distribution of the path of a simple symmetric random walk $(S_n)_{n=0,1,2,\dots}$ conditioned to stay positive, up to the last visit of $x_1$.
In the latter each such path is also possible, and they also happen with probability proportional to $2^{-k}$.
Now taking $Y_1=\inf\{n:S_{n+1}=S_n-1\}=\tilde S_{n'}$ with
$n'=1+\sup\{n:0\leqslant n < k, \tilde S_{n+1}= \tilde S_n+1\}$, it follows that
$a_1=Y_1$ whenever $\{n:S_{n+1}=S_n-1\}\cap
\llbracket 1,x_1-1 \rrbracket \ne\emptyset$.
The claim is proved.

As in the previous section, there is an i.i.d.\ sequence $(Y_1,Y_2,Y_3,\dots)$ with the property that the $n$-th step is successful if and only if the previous steps are successful and $a_{k-1}+Y_k\leqslant x_k$, in which case $a_k=a_{k-1}+Y_k$.

Now the expectation of $Y_1$ can be computed explicitly, and it is equal to $4$.
By the law of large numbers, the putative position $X_k=Y_1+\dots+Y_k$ of the $k$-th barrier grows like $4k$.
On the other hand, by the law of large numbers, the initial position $x_k$ of
the $k$-th particle to the right of the origin grows like $k/\mu$, and since
$\mu<\frac14$, with positive probability $X_k < x_k$ for all $k$.

Therefore, the construction is successful with positive probability and the proof is finished.
\end{proof}

\section{Concluding remarks}
\label{sec1concluding}

In our stabilizing strategy, the region to be avoided is the set of sites where
particles have been settled, together with the corrupted sites, or rather the
convex hull of this set.
In one dimension, the volume of this region is proportional to the number of
such particles.
This does not seem to have a direct analogue in higher dimensions:
straightforward application of a similar settling procedure would result in a
Diffusion-Limited Aggregation type of growth.
Thus, even though our approach seems promising for other settings, the proof of
local fixation is so far restricted to one dimension.

Yet, the construction presented in Section~\ref{sec1flag}, as well as Lemmas~\ref{lemma1wakeupmonotonicity},~\ref{lemma1halftopplingmonotonicity}, and Corollary~\ref{cor1abelianwithhalftopplings}, hold in any dimension.
In particular, $\mu_c$ is always well defined.
For the ARW, $\mu_c$ is non-decreasing in $\lambda$ by Lemma~\ref{lemma1wakeupmonotonicity}.

\bigskip

The main results have been stated and proved under simple conditions in order to
keep the presentation as transparent as possible.
Below we discuss some straightforward \emph{Generalizations}.

For the ARW, the proof that $\mu_c > 0$ when the jumps are to nearest neighbors can be adapted to bounded jumps without any difficulty.
In this case we cannot get the sharp estimate $\mu_c \geqslant \frac \lambda{1+\lambda}$, though we conjecture that it should still hold.
For the SSM with asymmetric jumps, the same proof still gives $\mu_c>0$, but the
lower bound for $\mu_c$ degenerates as $q \to 0$, where $q:=p(+1)=1-p(-1)$.

The Poissonian distribution $\nu$ of the initial amount of particles at a given site plays no special role in the proof of local fixation.
One can replace $\nu$ by any translation-invariant ergodic distribution with finite first moment $\mu$ and the proofs presented in the previous sections still give fixation as long as $\mu<\frac{\lambda}{1+\lambda}$ for the ARW or $\mu<\frac{1}{4}$ for the SSM.
For the existence of a $\mu_c\in[0,\infty]$ separating the absorbing and active phases we have used the stochastic ordering of the family $\{\nu_\mu\}_{0\leqslant\mu<\infty}$, parametrized by the density $\mu=\nu\big(\eta(0)\big)$.

The proof that $\mu_c \leqslant 1$ in one dimension generalizes with slight modifications to any finite-range random walk and any ergodic $\nu$.
Indeed, the proof shows that the system stays active for any ergodic $\nu$ with density $\mu>1$, or even $\mu=1$, as long as the fluctuations diverge, that is, $\limsup_N\sum_{x=0}^N|\eta(x)|-N=\infty$.
A proof that $\mu_c\leqslant 1$ in dimensions $d\geqslant2$ was recently given in~\cite{shellef-10}, using the framework described here and in~\cite{rolla-08}, and a more abstract proof was later given in~\cite{amir-gurelgurevich-10}.

\bigskip

We finish with some \emph{Open Problems}, see also~\cite{dickman-rolla-sidoravicius-10}.

It is remarkable that a proof of $\mu_c>0$ is lacking for any $d \geqslant 2$.
Whereas the trick of letting the particles evolve until hitting barriers and then setting traps worked well for $d=1$, for higher dimensions more involved steps will be necessary.

A proof that $\mu_c < 1$ for both the SSM and ARW remains as an open problem in any dimension.
For the ARW, it should hold for all $\lambda$, and moreover $\mu_c \to 0$ as $\lambda\to 0$.
Yet, even a proof that $\mu_c<1$ for \emph{some} $\lambda>0$ is missing.
For the SSM, mean-field theory and extensive simulations confirm the existence of sustained activity for densities greater than some critical value $\mu_c$, with $\mu_c< 1$ \cite{dickman-02}.
Simulations show that $\mu_c \approx 0.9489$~\cite{dickman-munoz-vespignani-zapperi-00}.

The only true parameter in the ARW model should be the density $\mu$, in the sense that there is a value $\mu_c$ (depending only on the dimension, halting rate, and jump probabilities) such that, for any ergodic initial distribution $\nu$, the system locally fixates if $\nu(\eta(0))<\mu_c$ and stays active if $\nu(\eta(0))>\mu_c$.

Besides universality of $\mu_c$ with respect to the distribution, a very
interesting open problem, but also very hard to approach, is whether there is
local fixation at $\mu=\mu_c$.
We believe that at the critical density $\mu=\mu_c$ the two models stay active (at least when the initial condition is i.i.d.\ with non-degenerate marginals), in marked contrast with several lattice models that exhibit phase transition, such as percolation or the Ising model, for which there is no percolation at criticality.
This has been proved only for the very particular case of totally asymmetric, nearest-neighbor walks on the one-dimensional lattice~\cite{hoffman-sidoravicius-04}.

\section*{Acknowledgments}

We thank R.~Dickman, who introduced to us these models, spent many hours of fruitful discussions and provided valuable numerical data.
We thank C.~Hoffman and A.~Teixeira for deep and inspiring discussions, and D.~Valesin for careful reading and valuable comments.
We also thank the anonymous referees whose comments considerably improved the
presentation.
L.~T.~Rolla thanks the hospitality of CWI, where part of this research was done.
This work had financial support from CNPq grant 141114/2004-5, FAPERJ, FAPESP grant 07/58470-1, and FSM-Paris.

\bibliographystyle{bib/siamrolla}
\bibliography{bib/leo}

\end{document}